\documentclass[leqno]{amsart}
\usepackage{amsfonts,amssymb,amsmath,amsgen,amsthm}
\usepackage{hyperref}

\theoremstyle{plain}
\newtheorem{theorem}{Theorem}[section]
\newtheorem{definition}[theorem]{Definition}

\newtheorem{lemma}[theorem]{Lemma}

\newtheorem{proposition}[theorem]{Proposition}

\theoremstyle{remark}
\newtheorem{remark}[theorem]{Remark}

\def\C{{\mathbf C}}
\def\R{{\mathbf R}}
\def\Sch{{\mathcal S}}
\def\O{\mathcal O}
\def\F{\mathcal F}

\def\({\left(}
\def\){\right)}
\def\<{\left\langle}
\def\>{\right\rangle}
\def\le{\leqslant}
\def\ge{\geqslant}
\newcommand{\chap}{\widehat}

\def\d{{\partial}}

\def\l{\lambda}

\def\si{{\sigma}}

\numberwithin{equation}{section}

\begin{document}

\title[Hartree equation in the Wiener algebra]{On the Cauchy problem
  for Hartree equation in the Wiener algebra}
\author[R. Carles]{R\'emi Carles}
\author[L. Mouzaoui]{Loun\`es Mouzaoui}
\address{CNRS \& Univ. Montpellier~2\\Math\'ematiques
\\CC~051\\34095 Montpellier\\ France}
\email{Remi.Carles@math.cnrs.fr}
\email{lounes.mouzaoui@univ-montp2.fr}
\begin{abstract}
We consider the mass-subcritical Hartree equation with a homogeneous
kernel, in the space of
square integrable functions whose Fourier transform is integrable. We
prove a global well-posedness result in this
space. On the other hand, we show that the Cauchy problem is not even
locally well-posed if we simply work in the space of functions whose
Fourier transform is integrable. Similar results are proven when the
kernel is not homogeneous, and is such that its Fourier transform
belongs to some Lebesgue space. 
\end{abstract}
\thanks{This work was supported by the French ANR project
  R.A.S. (ANR-08-JCJC-0124-01).}
\maketitle

\section{Introduction}
\label{sec:intro}
We consider the Cauchy problem for the following Hartree equation
\begin{equation}
  \label{eq:r3}
  i\d_t u + \Delta u =  \(K\ast |u|^2\)u,\quad t\in \R, \
  x\in\R^d;\quad u_{\mid t=0}=u_0,
\end{equation}
where $K$ denotes the Hartree kernel.
We  first deal with  the case of a  homogeneous kernel,
\begin{equation}
  \label{eq:noyau}
  K(x)=\frac{\l}{|x|^\gamma},\quad \l\in \R,\ \gamma>0.
\end{equation}
In \cite{Mou12}, it was proved that if $1\le d\le 3$ and $\gamma<d$,
the Cauchy problem \eqref{eq:r3} is locally well-posed in
$L^2(\R^d)\cap W$, where $W$ stands for the Wiener algebra (also
called Fourier algebra, according to the context)
\begin{equation*}
  W=\left\{f\in \Sch'(\R^d;\C),\quad \widehat f\in L^1(\R^d)\right\},
\end{equation*}
and the Fourier transform is defined, for $f\in L^1(\R^d)$, as
\begin{equation*}
  \widehat f(\xi) = \frac{1}{(2\pi)^{d/2}}\int_{\R^d}e^{-ix\cdot \xi}f(x)\mathrm{d}x.
\end{equation*}
In this paper, we investigate the global well-posedness for
\eqref{eq:r3}: we prove that if $d\ge 1$ and $\gamma<\min (2,d/2)$, then
the solution to \eqref{eq:r3}
is global in time in $L^2\cap W$. In view of the classical result
according to which \eqref{eq:r3} is globally well-posed in $L^2(\R^d)$,
our result can be understood as a propagation of the Wiener
regularity. On the other hand, the mere Wiener regularity does not
suffice to ensure even local well-posedness for \eqref{eq:r3}.

An advantage of working in $W$ lies in the fact that $W\hookrightarrow
L^\infty(\R^d)$, and, contrary to e.g. $H^s(\R^d$), $s>d/2$, $W$
scales like $L^\infty(\R^d)$ (note also that if $s>d/2$,
$H^s(\R^d)\hookrightarrow W$). So in a way, $W$ is the largest space
included in $L^\infty(\R^d)$ on which the Schr\"odinger group
$e^{it\Delta}$ acts continuously --- see Remark~\ref{rem:amalgam}
though. Recall that $e^{it\Delta}$ does not
map $L^\infty(\R^d)$ to itself, as shown by the explicit formula
$e^{i\Delta}(e^{-i|x|^2/4}) = \delta_{x=0}$, and the parabolic scale
invariance.
\begin{theorem}\label{theo:main}
  Let $d\ge 1$, $K$ given by \eqref{eq:noyau} with $\l\in \R$ and
  $0<\gamma<\min (2,d/2)$. If $u_0\in L^2(\R^d)\cap W$,  then \eqref{eq:r3}
  has a unique, global solution $u\in C(\R; L^2\cap W)$. In addition,
  its $L^2$-norm is conserved,
  \begin{equation*}
    \|u(t)\|_{L^2}=\|u_0\|_{L^2},\quad \forall t\in \R.
  \end{equation*}
\end{theorem}
\begin{remark}\label{rem:amalgam}
  One might be tempted to consider a space included in
$L^\infty$ which also
scales like $L^\infty$,  and which is
larger than $W=\F L^1$, namely the amalgam space $W(\F L^1;L^\infty)$. This
space consists essentially of functions which are locally in $W$,
and globally in $L^\infty$ (see e.g. \cite{CoNi08,CoZu10} for a
precise definition). Strichartz estimates in amalgam spaces have been
established in \cite{CoNi08} (even though the case $W(\F L^1;L^\infty)$
can never be considered). However, since the map $x\mapsto e^{i|x|^2}$
belongs to
$W(\F L^1;L^\infty)$, we see that $e^{it\Delta}$ does not act
continuously on
$W(\F L^1;L^\infty)$.
\end{remark}

We next show that \eqref{eq:r3} is not well-posed
in the mere Wiener algebra. Precisely have the following theorem:
\begin{theorem}
\label{th:ill}
Let $d\ge 1$, $K$ given by \eqref{eq:noyau} with $0<\gamma<d$. The Cauchy problem \eqref{eq:r3} is
locally well-posed in $W\cap L^2$, but not in $W$:
for all ball $B$ of $W$, for all $T>0$ the solution map
$ \varphi \in B \mapsto \psi \in C([0,T]; W)$ is not uniformly continuous.
\end{theorem}
\begin{remark}
In the case of the nonlinear Schr\"odinger equation
\begin{equation}\label{eq:nls}
i\partial_{t}u+\Delta u = \lambda|u|^{2\si}u,
\end{equation}
where $\si $ is an integer and $\lambda \in \R$, the Cauchy problem is
locally well-posed in $W$ (see \cite{CDS10}). From the above result,
this is in sharp contrast with the case of the Hartree equation. On
the other hand, it is not clear that the Cauchy problem for
\eqref{eq:nls} is \emph{globally} well-posed in $L^2(\R^d)\cap W$,
even in the case $d=\si=1$. We note that in \cite{HyTs12}, the authors
study \eqref{eq:nls} in the one-dimensional case $d=1$, with $\si<2$
not necessarily an integer. They prove local well-posedness in
\begin{equation*}
  \hat L^p= \{f\ : \ \hat f\in L^{p'}\},
\end{equation*}
for $p$ in some open neighborhood of $2$. Global well-posedness
results for initial data in $\hat L^p$ are established, in spaces
based on dispersive estimates.
\end{remark}
The kernel $K$ given by
\eqref{eq:noyau} is such that its Fourier transform belongs to no
Lebesgue space, but to a weak Lebesgue space, from the following property (see
e.g. \cite[Proposition~1.29]{BCD11}):
\begin{proposition}\label{prop:hatK}
  Let $d\ge 1$ and $0<\gamma<d$. There exists $C=C(\gamma,d)$ such
  that the Fourier transform of $K$ defined by \eqref{eq:noyau} is
  \begin{equation*}
    \widehat K(\xi) =  \frac{\lambda C}{|\xi|^{d-\gamma}}.
  \end{equation*}
\end{proposition}
The final result of this paper is concerned with the case where the
kernel $K$ is such that its Fourier transform belongs to some Lebesgue space.
\begin{theorem}\label{theo:Lp}
  Let $d\ge 1$.
  \begin{itemize}
  \item Let $p\in [1,\infty]$, and suppose that $K$ is such that $\widehat
    K\in L^p(\R^d)$.  If
  $u_0\in L^2(\R^d)\cap W$,  then there exists $T>0$ such that
    \eqref{eq:r3}  has a unique solution $u\in C([-T,T]; L^2\cap W)$.
In addition, its $L^2$-norm is conserved, 
  \begin{equation*}
    \|u(t)\|_{L^2}=\|u_0\|_{L^2},\quad \forall t\in [-T,T].
  \end{equation*}
If $K\in W$ ($p=1$), then the solution is global: $u\in C(\R; L^2\cap W)$.
\item Suppose that $K$ is such that $\widehat K\in L^\infty(\R^d)$. If
  $u_0\in W$,  then there exists $T>0$ such that
    \eqref{eq:r3}  has a unique solution $u\in C([-T,T]; W)$.
\item For any $p\in [1,\infty)$, one can find $K$ with  $\widehat K\in
  L^p(\R^d)\setminus L^\infty(\R^d)$, such that
  \eqref{eq:r3} is locally well-posed in $L^2(\R^d)\cap W$, but not
  in $W$: for all ball $B$ of 
  $W$, for all $T>0$ the solution map 
$ \varphi \in B \mapsto \psi \in C([0,T]; W)$ is not uniformly continuous.
  \end{itemize}
\end{theorem}

\section{Standard existence results and properties}
\label{sec:standard}

\subsection{Main properties of the Wiener algebra}
\label{sec:wiener}

The space  $W$ enjoys the following elementary properties (see
\cite{CoLa09,CDS10}):
\begin{enumerate}
\item $W$ is a Banach space, continuously embedded into
$L^\infty(\R^d)$.
\item $W$ is an algebra, and the
  mapping $(f,g)\mapsto fg$ is continuous from $W^2$ to
  $W$, with
  \begin{equation*}
    \|fg\|_{W}\le \|f\|_{W}\|g\|_{W},\quad \forall f,g\in W.
  \end{equation*}
\item  For all $t\in \R$, the free Schr\"odinger group
  $e^{i t \Delta}$
is unitary on $W$.
\end{enumerate}

\subsection{Existence results based on Strichartz estimates}
\label{sec:stri}

For the sake of completeness, we recall standard definition and
results.
\begin{definition}\label{def:adm}
 A pair $(p,q)\not =(2,\infty)$ is admissible if $p\ge 2$, $q\ge  2$,
 and
$$\frac{2}{p}= d\left( \frac{1}{2}-\frac{1}{q}\right).$$
\end{definition}

\begin{proposition}[From \cite{GV85c,KT}]\label{prop:strichartz}
  $(1)$ For any admissible pair $(p,q)$, there exists $C_{q}$  such that
$$
\|e^{it\Delta} \varphi\|_{L^{p}(\R;L^{q})} \le C_q
\|\varphi \|_{L^2},\quad \forall \varphi\in L^2(\R^d).
$$
$(2)$ Denote
\begin{equation*}
  D(F)(t,x) = \int_0^t e^{i(t-\tau)\Delta}F(\tau,x)\mathrm{d}\tau.
\end{equation*}
For all admissible pairs $(p_1,q_1)$ and~$
    (p_2,q_2)$,  there exists $C=C_{q_1,q_2}$
    such that for all interval $I\ni 0$,
\begin{equation*}
      \left\lVert D(F)
      \right\rVert_{L^{p_1}(I;L^{q_1})}\le C \left\lVert
      F\right\rVert_{L^{p'_2}\(I;L^{q'_2}\)},\quad \forall F\in
    L^{p'_2}(I;L^{q'_2}).
\end{equation*}
\end{proposition}

\begin{proposition}\label{prop:r3L2}
  Let $d\ge 1$, $K$ given by
  \eqref{eq:noyau} with $\l\in \R$ and $0<\gamma<\min(2,d)$. If
  $u_0\in L^2(\R^d)$, then \eqref{eq:r3} has a 
  unique, global solution
  \begin{equation*}
    u\in C(\R;L^2)\cap L^{8/\gamma}_{\rm loc}(\R;L^{4d/(2d-\gamma)} ).
  \end{equation*}
In addition,
  its $L^2$-norm is conserved,
  \begin{equation*}
    \|u(t)\|_{L^2}=\|u_0\|_{L^2},\quad \forall t\in \R,
  \end{equation*}
and for all admissible pair $(p,q)$, $u \in L^p_{\rm
  loc}(\R;L^q(\R^d))$.
\end{proposition}
\begin{proof}
We give the main technical steps of the proof,  and refer to  \cite{CazCourant}
for details.
 By Duhamel's formula, we write \eqref{eq:r3} as
\begin{equation*}
u(t)=e^{it\Delta}u_0-
i\int_0^{t}e^{i(t-\tau)\Delta}(K\ast |u|^{2}u)(\tau)\mathrm{d}\tau=:\Phi(u)(t).
\end{equation*}
Introduce the space
 \begin{equation*}
\begin{aligned}
Y(T)&=\{\phi\in C([0,T]; L^{2}(\R^{d})):
\|\phi\|_{L^{\infty}([0,T];\ L^{2}(\R^{d}))}\le
2\|u_0\|_{L^{2}(\R^{d})},\\
&\quad \|\phi\|_{L^{8/\gamma}([0,T];\
L^{4d/(2d-\gamma)}(\R^{d}))}\le 2C(8/\gamma)\|u_0\|_{L^{2}(\R^{d})}\},
\end{aligned}
\end{equation*}
and the distance $$ d(\phi_{1},
\phi_{2})=\|\phi_{1}-\phi_{2}\|_{L^{8/\gamma}([0,T];\
L^{4d/(2d-\gamma)})},$$
where $C(8/\gamma)$ stems from Proposition~\ref{prop:strichartz}.
Then $(Y(T), d)$ is a complete metric space,
  as remarked in \cite{Kato87} (see also \cite{CazCourant}).
Hereafter, we denote by
\begin{equation*}
q=\frac{8}{\gamma}, \quad r=\frac{4d}{2d-\gamma}, \quad
\theta=\frac{8}{4-\gamma},
\end{equation*}
and $\|\cdot\|_{_{L^{a}([0,T];\ L^{b}(\R^{d}))}}$ by
$\|\cdot\|_{L^{a}L^{b}}$ for simplicity. Notice that $\left(q, r\right)$ is
admissible and
\begin{equation*}
\frac{1}{q^{\prime}}=\frac{4-\gamma}{4}+\frac{1}{q}=\frac{1}{2}+\frac{1}{\theta}
\quad ; \quad
\frac{1}{r^{\prime}}=\frac{\gamma}{2d}+\frac{1}{r}\quad;\quad
\frac{1}{2}=\frac{1}{\theta}+\frac{1}{q}.
\end{equation*}
 By using Strichartz estimates, H\"older inequality and
 Hardy--Littlewood--Sobolev inequality, we have, for $(\underline
 q,\underline r)\in \{(q,r),(\infty,2)\}$:
\begin{equation*}
\begin{aligned}
\|\Phi(u)\|_{L^{\underline q}L^{\underline r}}&\le C(\underline q)
\|u_0\|_{L^{2}}+C(\underline q,q)\left\lVert (K\ast|u|^{2})u
\right\rVert_{L^{q'}L^{r'}}\\
&\le C(\underline q)\|u_0\|_{L^{2}}+
C(\underline q,q)\left\lVert
 K\ast |u|^{2}\right\rVert_{L^{4/(4-\gamma)}L^{2d/\gamma}}
\|u\|_{L^{q}L^{r}}\\
& \le C(\underline q)
\|u_0\|_{L^{2}}+C\|u\|_{L^{\theta}L^{r}}^{2}\|u\|_{L^{q}L^{r}}\\
&\le C(\underline q)\|u_0\|_{L^{2}}+CT^{1-\gamma/2}\|u\|_{L^{q}L^{r}}^{3},
\end{aligned}
\end{equation*}
 for any $u\in Y(T)$, with $C(\infty)=1$ by the standard energy
 estimate. To show the contraction property
of $\Phi$,  for any $v, w\in
Y(T)$, we get
\begin{align*}
\|\Phi(v)-\Phi(w)\|_{L^{q}L^{r}}
&\lesssim\|K\ast |v|^{2}\|_{L^{4/(4-\gamma)}L^{2d/\gamma}}
\|v-w\|_{L^{q}L^{r}}\\
&\quad
+\left\lVert K\ast \left\lvert |v|^{2}-|w|^{2}\right\rvert
\right\rVert_{L^{2}L^{2d/\gamma}}\|w\|_{L^{\theta}L^{r}}\\
&\lesssim
\left(\|v\|_{L^{\theta}L^{r}}^{2}+\|w\|_{L^{\theta}L^{r}}^{2}\right)\|v-w\|_{L^{q}L^{r}}\\
& \le
C T^{1-\gamma/2}(\|v\|_{L^{q}L^{r}}^{2}+\|w\|_{L^{q}L^{r}}^{2})\|v-w\|_{L^{q}L^{r}}.
\end{align*}
 Thus $\Phi$ is a contraction from $Y(T)$ to $Y(T)$ provided
that $T$ is sufficiently small. Then there exists  a unique $u\in
Y(T)$ solving \eqref{eq:r3}. The global existence of the
solution for \eqref{eq:r3} follows from the conservation
of the $L^{2}$-norm of $u$. The last property of the proposition then
follows from Strichartz estimates applied with an arbitrary admissible
pair on the left hand side, and the same pairs as above on the right
hand side.
\end{proof}

\section{Proof of Theorem~\ref{theo:main}}
\label{sec:inter}

Thoughout this section, we assume that the kernel $K$ is given by
\eqref{eq:noyau}. 
\subsection{Uniqueness}
\label{sec:unique}
Uniqueness stems from the local well-posedness result established in
\cite{Mou12}, based on the following lemma, whose proof is recalled
for the sake of completeness.

\begin{lemma}\label{lem:estW}
  Let $0<\gamma<d$. There exists $C$ such that for all $f,g\in
  L^2(\R^d)\cap W$,
  \begin{equation*}
    \left\| \(K\ast |f|^2\)f -\(K\ast |g|^2\)g\right\|_{L^2\cap W}
    \le C\(\|f\|_{L^2\cap W}^2 + \|g\|_{L^2\cap W}^2\)\|f-g\|_{L^2\cap W}.
  \end{equation*}
\end{lemma}
\begin{proof}
  Let $\kappa_1= {\bf 1}_{\{|\xi|\le 1\}}\widehat K$ and $\kappa_2=
  {\bf 1}_{\{|\xi|>1\}}\widehat K$. In view of
  Proposition~\ref{prop:hatK}, $\kappa_1\in L^p(\R^d)$ for all $p\in
  [1,\frac{d}{d-\gamma})$ and $\kappa_2\in L^q(\R^d)$ for all $q\in
  (\frac{d}{d-\gamma},\infty]$. For $h\in L^1(\R^d)\cap W$, we have
  \begin{align*}
    \|K\ast h\|_{W} & \lesssim \|\kappa_1 \widehat h\|_{L^1} + \|\kappa_2
    \widehat h\|_{L^1} \lesssim \|\kappa_1\|_{L^1}\|\widehat
    h\|_{L^\infty} +\|\kappa_2\|_{L^\infty}\|\widehat    h\|_{L^1}\\
   & \lesssim \|h\|_{L^1} + \|\widehat    h\|_{L^1},
  \end{align*}
where we have used Hausdorff--Young inequality. Writing
\begin{align*}
  \(K\ast |f|^2\)f -\(K\ast |g|^2\)g = \(K\ast |f|^2\)(f-g) +\(K\ast
  (|f|^2-|g|^2)\)g,
\end{align*}
the lemma follows, since $W$ is a Banach algebra embedded
into $L^\infty(\R^d)$.
\end{proof}
We infer uniqueness in $L^2(\R^d)\cap W$ for \eqref{eq:r3} as soon as
$0<\gamma<d$:
\begin{proposition}\label{prop:uniqueness}
Let $0<\gamma<d$, $T>0$, and  $u,v\in C([0,T];L^2\cap W)$ solve \eqref{eq:r3},
with the same initial datum $u_0\in L^2(\R^d)\cap W$. Then $u\equiv v$.
\end{proposition}
\begin{proof}
 Duhamel's
formula yields
\begin{equation*}
  u(t)-v(t) = -i \int_0^t e^{i(t-\tau)\Delta}\(\(K\ast |u|^2\)u
  -\(K\ast |v|^2\)v\)(\tau)\mathrm{d}\tau.
\end{equation*}
Since the Schr\"odinger group is unitary on $L^2$ and on $W$, Minkowski
inequality and Lemma~\ref{lem:estW} yield, for $t\ge 0$,
\begin{align*}
  \|u(t)-v(t)&\|_{L^2\cap W} \lesssim \int_0^t \(\|u(\tau)\|_{L^2\cap
    W}^2 +\|v(\tau)\|_{L^2\cap    W}^2 \)\|u(\tau)-v(\tau)\|_{L^2\cap
    W}\mathrm{d}\tau\\
&\lesssim \( \|u\|_{L^\infty([0,T];L^2\cap W)}^2
+\|v\|_{L^\infty([0,T];L^2\cap W)}^2\) \int_0^t \|u(\tau)-v(\tau)\|_{L^2\cap
    W}\mathrm{d}\tau.
\end{align*}
Gronwall lemma implies $u\equiv v$.
\end{proof}
\subsection{Existence}
\label{sec:existence}

In view of Lemma~\ref{lem:estW}, the standard fixed point argument
yields:

\begin{proposition}\label{prop:local}
  Let $d\ge 1$, $\l\in \R$, $0<\gamma<d$,  and $K$ given by
  \eqref{eq:noyau}. If $u_0\in L^2(\R^d)\cap W$, then there exists
  $T>0$ depending only on $\lambda,\gamma,d$ and $\|u_0\|_{L^2\cap
    W}$, and a unique $u\in C([0,T];L^2\cap W)$ to \eqref{eq:r3}.
\end{proposition}

Taking Proposition~\ref{prop:r3L2} into account, to establish
Theorem~\ref{theo:main},  it suffices to prove that the Wiener norm of
$u$ cannot become unbounded in finite time.

Resuming the decomposition of $\widehat K$ introduced in the proof of
Lemma~\ref{lem:estW}, we find
\begin{align*}
  \|u(t)\|_{W} &\le \|u_0\|_{W} + \int_0^t \left\| \(K\ast
    |u(\tau)|^2\)u(\tau)\right\|_{W} \mathrm{d}\tau \\
&\le \|u_0\|_{W} + \int_0^t \left\| K\ast
    |u(\tau)|^2\right\|_{W}\|u(\tau)\|_{W} \mathrm{d}\tau \\
&\le \|u_0\|_{W} + \int_0^t \(\| \kappa_1\|_{L^1}
    \|u(\tau)\|_{L^2}^2+
    \|\kappa_2\|_{L^p}\left\|\widehat{|u(\tau)|^2}\right\|_{L^{p'}}\)
    \|u(\tau)\|_{W} \mathrm{d}\tau ,
\end{align*}
provided that $p>\frac{d}{d-\gamma}$. Using the conservation of the
$L^2$-norm of $u$ and Hausdorff--Young inequality, we infer, if $p\le
2$:
\begin{align*}
 \|u(t)\|_{W} \lesssim  \|u_0\|_{W} +  \int_0^t \(\| \kappa_1\|_{L^1}
    \|u_0\|_{L^2}^2+
    \|\kappa_2\|_{L^p}\left\||u(\tau)|^2\right\|_{L^{p}}\)
    \|u(\tau)\|_{W} \mathrm{d}\tau.
\end{align*}
To summarize, for all $1<\frac{d}{d-\gamma}<p\le 2$, there exists $C$
such that
\begin{equation*}
  \|u(t)\|_{W} \le  \|u_0\|_{W} + C \int_0^t \|u(\tau)\|_{W} \mathrm{d}\tau +
  C\int_0^t \left\|u(\tau)\right\|^2_{L^{2p}}
    \|u(\tau)\|_{W} \mathrm{d}\tau.
\end{equation*}
The above requirement on $p$ can be fulfilled if and only if
$0<\gamma<d/2$.
To take advantage of Proposition~\ref{prop:r3L2}, introduce $\alpha>
1$
such that $(2\alpha,2p)$ is admissible. This is possible provided that
$2p<\frac{2d}{d-2}$ when $d\ge 3$: this condition is compatible with
the requirement $p>\frac{d}{d-\gamma}$ if and only if
$\gamma<2$. Using H\"older inequality for the
last integral, we have
\begin{equation*}
  \|u(t)\|_{W} \le  \|u_0\|_{W} + C \int_0^t \|u(\tau)\|_{W} \mathrm{d}\tau +
  C\left\|u\right\|^2_{L^{2\alpha}([0,t];L^{2p})}
    \|u\|_{L^{\alpha'}([0,t];W)}.
\end{equation*}
Set
\begin{equation*}
  \omega (t) =\sup_{0\le \tau\le t}\|u(\tau)\|_{W}.
\end{equation*}
For a given $T>0$, $\omega$ satisfies an
estimate of the form
\begin{equation*}
  \omega(t) \le \|u_0\|_{W} + C\int_0^t\omega(\tau)\mathrm{d}\tau +
  C_0(T)\(\int_0^t \omega(\tau)^{\alpha'}\mathrm{d}\tau\)^{1/\alpha'},
\end{equation*}
provided that $0\le t\le T$, and where we have used the fact that
$\alpha'$ is finite.
Using H\"older inequality, we infer
\begin{equation*}
  \omega(t) \le \|u_0\|_{W} +
  C_1(T)\(\int_0^t \omega(\tau)^{\alpha'}\mathrm{d}\tau\)^{1/\alpha'}.
\end{equation*}
Raising the above estimate to the power $\alpha'$, we find
\begin{equation*}
  \omega(t)^{\alpha'} \lesssim 1 + \int_0^t \omega(\tau)^{\alpha'}\mathrm{d}\tau.
\end{equation*}
Gronwall lemma shows that $\omega \in L^\infty([0,T])$. Since $T>0$ is
arbitrary, $\omega\in L^\infty_{\rm loc}(\R)$, and the result
follows.

\section{Ill-posedness in the mere Wiener algebra}
\label{sec:ill}
In this section we still assume that $K$ is given by
\eqref{eq:noyau}. We show that the Cauchy problem \eqref{eq:r3} is ill-posed
in the mere Wiener algebra, i.e without including $ L^{2}$.
We recall the definition of well-posedness for the problem \eqref{eq:r3}.
\begin{definition}
\label{well}
Let $(S, \|\cdot\|_{S})$ be a Banach space of initial data, and
$(D,\|\cdot\|_{D})$ be a Banach space of space-time functions.
The Cauchy problem \eqref{eq:r3} is well posed from $D$ to $S$ if, for all bounded
subset $B\subset D$, there exist $T > 0$ and a Banach
space $X_{T} \hookrightarrow  C([0,T]; S)$ such that:
\begin{enumerate}
\item For all $\varphi\in B$, \eqref{eq:r3} has a unique solution $ \psi\in X_{T}$
with $\psi|_{t=0}=\varphi$.
\item The mapping $B\ni\varphi\mapsto \psi\in C([0,T];S)$ is uniformly
continuous.
\end{enumerate}
\end{definition}

\begin{proof}[Proof of Theorem~\ref{th:ill}]
In view of \cite{BeTa06}, it suffices to prove that one term in the Picard
iterations of $\Phi$ defined 
in Section~\ref{sec:standard} does not verify the Definition~\ref{well}.
We argue by contradiction and assume that \eqref{eq:r3} is well-posed from
$W$ to $W$. Then, let $T>0$ be the local time existence of the solution. We recall that for $0<\gamma<d$, \eqref{eq:r3} is
well-posed from $W\cap L^{2}$ to $W\cap L^{2}$ (see
\cite[Theorem~2.1]{Mou12}, recalled in Proposition~\ref{prop:local}).
We define the following operator associated to the second Picard iterate:
\[
D(f)(t,x)=-i\int_{0}^{t} e^{i(t-\tau)\Delta}(K*|\psi|^{2}\psi)(\tau,x)\; \mathrm{d}x,
\]
where $\psi$ is the solution of the free equation
\[
i\partial_{t}\psi + \Delta \psi = 0; \quad \psi |_{t=0}=f,
\]
that is $\psi(t) = e^{it\Delta}f$. We denote $e^{it\Delta}f=U(t)f$.
By \cite[Proposition 1]{BeTa06} the operator $D$ is continuous
from $W$ to $W$, that is,
\begin{equation}
\label{ill-control}
\|D(f)(t)\|_{W}\le C\|f\|^{3}_{W},\quad \forall t\in [0,T],
\end{equation}
for some positive constant $C$. Let $f\in \Sch(\R^{d})$ be an element of
the Schwartz space. We define a family of functions indexed by $h >0$
by
\[
f^{h}(x)=f(hx).
\]
For all $h>0, \|f^{h}\|_{L^{2}}=\frac{1}{h^{d/2}}\|f\|_{L^{2}}$ so for $h>0$ close to $0$
the family $(f^{h})_{h>0}$ leaves any compact of $L^{2}$. Remark that
$\chap{f^{h}}(\xi)=h^{-d}\chap{f}(\frac{\xi}{h})$ and
$ \|f^{h}\|_{W}=\|f\|_{W} < \infty$, so \eqref{ill-control} yields
\begin{equation}
\label{fh_W}
\|D(f^{h})(t)\|_{W}\le C\|f\|^{3}_{W}.
\end{equation}
We develop the expression of $\|D(f^{h})(t)\|_{W}$. We have
\[
\begin{aligned}
&\mathcal{F}(K*|U(\tau)f^{h}|^{2}U(\tau)f^{h})(\xi)=\phantom{\int_{\R^{d}}
\int_{\R^{d}}\chap{K}(\xi-y) e^{i\tau|\xi-y-z|^{2}}\chap{f^{h}}(\xi-y-z)
e^{i\tau|y|^{2}}\chap{f^{h}}(y)\;\mathrm{d}y}\\
&=(2\pi)^{d/2}\mathcal{F}(K*|U(\tau)f^{h}|^{2})*\mathcal{F}(U(\tau)f^{h})(\xi)\\
&=(2\pi)^{d/2}\int_{\R^{d}} \mathcal{F}(K*|U(\tau)f^{h}|^{2})(\xi-y)
\mathcal{F}(U(\tau)f^{h})(y)\;\mathrm{d}y\\
&=\int_{\R^{d}} e^{i\tau|y|^{2}}\chap{K}(\xi-y)\mathcal{F}
(|U(\tau)f^{h}|^{2})(\xi-y)\chap{f^{h}}(y)\;\mathrm{d}y\\
&=(2\pi)^{d/2}\iint_{\R^{2d}} e^{i\tau|y|^{2}}
e^{i\tau|\xi-y-z|^{2}}e^{-i\tau|z|^{2}}\chap{K}(\xi-y)
\chap{f^{h}}(\xi-y-z)\chap{\overline{f^{h}}}(z)
\chap{f^{h}}(y)\;\mathrm{d}y\mathrm{d}z\\
&=\frac{(2\pi)^{d/2}}{h^{3d}}\iint_{\R^{2d}}
e^{i\tau|y|^{2}}
e^{i\tau|\xi-y-z|^{2}}e^{-i\tau|z|^{2}}\chap{K}(\xi-y)
\chap{f}\(\frac{\xi-y-z}{h}\)\chap{\overline{f}}\(\frac{z}{h}\)
\chap{f}\(\frac{y}{h}\)\;\mathrm{d}y\mathrm{d}z.\\
\end{aligned}
\]
Taking the $W$-norm gives
\[
\begin{aligned}
\|D(f^{h})(t)\|_{W}&=\|\chap{D(f^{h})}(t)\|_{L^{1}}\\
&=\int_{\R^{d}}\left|\int_{0}^{t} \mathcal{F}(U(t-\tau)(K*|U(\tau)f^{h}|^{2}
U(\tau)f^{h})(\xi)\;\mathrm{d}\tau\right|\;\mathrm{d}\xi\\
&=\int_{\R^{d}}\left|\int_{0}^{t} e^{i(t-\tau)|\xi|^{2}}
\mathcal{F}(K*|U(\tau)f^{h}|^{2}U(\tau)f^{h})(\xi)\;\mathrm{d}\tau\right|
\;\mathrm{d}\xi.\\
\end{aligned}
\]
We replace $\mathcal{F}(K*|U(\tau)f^{h}|^{2}U(\tau)f^{h})(\xi)$ by
its integral
formula above and apply the following changes of variable:
$\xi'=\xi/h, y'=y/h, z'=z/h, \tau'=\tau h^{2}$. We obtain
\begin{equation}
\label{Dfh_W}
\|D(f^{h})(t)\|_{W}=\frac{1}{h^{d-\gamma+2}}\|D(f)(th^{2})\|_{W}.
\end{equation}
Let $s\in (0, T)$. We examine more closely the term $F(s):=\int_{0}^{s}
U(s-\tau)g(\tau)\mathrm{d}\tau$ where $g(s):= (K*|U(s)f|^{2})U(s)f.$
Taylor formula  yields
\[
F(s)=F(0)+F'(0)s + \frac{s^{2}}{2}\int_{0}^{1} (1-\theta)F''(s\theta)\;\mathrm{d}\theta.
\]
We have $F(0)=0$, so for $s\in [0,1]$,
\[
\begin{aligned}
\|F(s)-F'(0)s\|_{W} &\le s^{2}\left\|\int_{0}^{1}
  (1-\theta)F''(s\theta)\;\mathrm{d}\theta \right\|_{W}\;\mathrm{d}\theta\\
& \le s^{2} \int_{0}^{1} \|F''(s\theta)\|_{W} \;\mathrm{d}\theta
\le s^{2} \|F''\|_{L^{\infty}([0,1];W)}.
\end{aligned}
\]
The first and second derivatives of $F$ are given by
\begin{align*}
F'(s)&=g(s)+\int_{0}^{s}U(s-\tau)i\Delta g(\tau)\;\mathrm{d}\tau,\\
F''(s)&=g'(s)+ i\Delta g(s) - \int_{0}^{s}
U(s-\tau)\Delta^{2}g(\tau)\;\mathrm{d}\tau,
\end{align*}
so $F'(0)=g(0)$ and
\[
\|F''\|_{L^{\infty}([0,1];W)} \le \|g'\|_{L^{\infty}([0,1];W)} +
\|\Delta g\|_{L^{\infty}([0,1];W)}+
\|\Delta^{2}g\|_{L^{\infty}([0,1];W)}.
\]
From the formula of $g$, since $f \in \mathcal{S}$ we can see easily
that $F''\in L^{\infty}([0,1];W)$ (uniformly in $h\in (0,1]$).
We obtain
\[
\|F(s)-sg(0)\|_{W} \le C s^{2},
\]
where $C$ is a positive constant independent of $s$ (it depends on $f$
and $\gamma$). Thus,
\[
\|F(s)\|_{W}= s\|g(0)\|_{W} + \O(s^{2}).
\]
In particular, for all $t, h>0$
\[
\|D(f)(th^{2})\|_{W}=\|F(th^{2})\|_{W}= th^{2}\|g(0)\|_{W} + \O\(t^{2}h^{4}\).
\]
This implies that
\[
\|D(f^{h})(t)\|_{W}=\frac{1}{h^{d-\gamma+2}}\|D(f)(th^{2})\|_{W}=
\frac{t}{h^{d-\gamma}}\|g(0)\|_{W}+\frac{1}{h^{d-\gamma}}\O\(t^{2}h^{2}\).
\]
Fix $t>0$:
\[
\lim\limits_{h\to 0^{+}}
\frac{t}{h^{d-\gamma}}\|g(0)\|_{W}+\frac{1}{h^{d-\gamma}}\O\(t^{2}h^{2}\)= +\infty.
\]
We deduce that for $h>0$ sufficiently close to $0$,
\[
\|D(f^{h})(t)\|_{W} > C\|f\|_{W}^{3}.
\]
This contradicts
\eqref{ill-control}, and Theorem \ref{th:ill} follows.
\end{proof}

\section{Proof of Theorem~\ref{theo:Lp}}
 We decompose the proof of Theorem~\ref{theo:Lp}
 according to the three cases considered.

\subsection{Well-posedness in $L^{2}\cap W$}
We assume that
$K$ is such that $\chap{K}\in L^{p}$ for some $p\in [1,\infty]$, and
we consider an initial
data $u_{0}\in L^{2}\cap W$. 
For $T > 0$ we define the following space
\[
E_{T}=\{u\in C([0,T];L^{2}\cap W),\quad \|u\|_{L^{\infty}([0,T];L^{2}\cap W)}\le 2\|u_{0}\|_{L^{2}\cap W}\}.
\]
It is a complete space metric when equipped with the metric
\[
d(u,v)= \|u-v\|_{L^{\infty}([0,T];L^{2}\cap W)}.
\]

\begin{lemma}\label{estWL^2}
Let $p\in [1,\infty]$, and $K$ such that $\widehat K\in L^p$. There exists
$C$ such that for all $f, g\in L^{2}\cap W$, 
\[
\|K*(fg)\|_{W} \le C\|\chap{K}\|_{L^{p}} \|f\|_{L^{2}\cap W}\|g\|_{L^{2}\cap W}.
\]
\end{lemma}

\begin{proof}
Let $ f, g\in L^{2}\cap W$. We denote by $p'$ the H\"older conjugate
exponent of $p$. We have
\[
\begin{aligned}
\|K*(fg)\|_{W} = (2\pi)^{d/2} \|\chap{K}\; \chap{fg}\|_{L^{1}}
&\le C\|\chap{K}\|_{L^{p}}\|\chap{fg}\|_{L^{p'}}
\le C\|\chap{K}\|_{L^{p}} \|\chap{f}*\chap{g}\|_{L^{p'}}\\
&\le C \|\chap{K}\|_{L^{p}} \|\chap{f}\|_{L^{q}}\|\chap{g}\|_{L^{q}},
\end{aligned}
\]
where $1+\frac{1}{p'}= \frac{2}{q}$. Since $q\in [1,2]$ there exists
$\theta \in [0,1]$ such that $\frac{1}{2}+ \frac{\theta}{2}= \frac{1}{q}$
and
\[
\|\chap{f}\|_{L^{q}} \le \|\chap{f}\|^{\theta}_{L^{1}}\|\chap{f}\|^{1-\theta}_{L^{2}}
=\|f\|^{\theta}_{W}\|f\|^{1-\theta}_{L^{2}}
\le \|f\|_{L^{2}\cap W}.
\]
The lemma follows.
\end{proof}

Let $\Phi$ as defined in Section~\ref{sec:standard}. In view of the
previous lemma we have, for $u\in E_{T}$ and $t\in [0,T]$, 
\[
\begin{aligned}
\|\Phi(u)(t)\|_{L^{2}\cap W} &\le \|u_{0}\|_{L^{2}\cap W} +
\int_{0}^{t} \|K*|u(\tau)|^{2}u(\tau)\|_{L^{2}\cap W}\;\mathrm{d}\tau\\
&\le \|u_{0}\|_{L^{2}\cap W} +
\int_{0}^{t} \|K*|u(\tau)|^{2}\|_{W}\|u(\tau)\|_{L^{2}\cap W}\;\mathrm{d}\tau\\
&\le \|u_{0}\|_{L^{2}\cap W} +
C\|\chap{K}\|_{L^{p}}\int_{0}^{t} \|u(\tau)\|^{3}_{L^{2}\cap W} \;\mathrm{d}\tau.\\
\end{aligned}
\]
We obtain
\[
\|\Phi(u)(t)\|_{L^{\infty}([0,T];L^{2}\cap W)} \le \|u_{0}\|_{L^{2}\cap W}
+ C\|\chap{K}\|_{L^{p}}\|u_{0}\|^{3}_{L^{\infty}([0,T];L^{2}\cap W)}T.
\]
For $T$ sufficiently small (depending on $\|u_{0}\|_{L^{2}\cap W}$)
$\|\Phi(u)\|_{L^{\infty}_{T}L^{2}} \le 2\|u_{0}\|_{L^{2}\cap W}$. 
Let $u,v \in E_{T}$. From Lemma \ref{estWL^2} we have
\[
\begin{aligned}
&\|\Phi(u)(t)-\Phi(v)(t)\|_{L^{2}\cap W}
\\
&\le \int_{0}^{t}
\|K*(|u(\tau)|^{2}-|v(\tau)|^{2})u(\tau)\|_{L^{2}\cap
  W}\;\mathrm{d}\tau\\
&\quad+\int_{0}^{t}\|K*|v(\tau)|^{2}(u(\tau)-v(\tau)\|_{L^{2}\cap
  W}\;\mathrm{d}\tau\\ 
&\le \int_{0}^{t} \|K*(|u(\tau)|^{2}-|v(\tau)|^{2})\|_{W}\|u(\tau)\|_{L^{2}\cap W}\;\mathrm{d}\tau\\
&\quad+ \int_{0}^{t} \|K*|v(\tau)|^{2}\|_{W}\|u(\tau)-v(\tau)\|_{L^{2}\cap W}\;\mathrm{d}\tau\\
&\le\int_{0}^{t} C\|\chap{K}\|_{L^{p}} \(\|u(\tau)\|^{2}_{L^{2}\cap
  W}+
\|v(\tau)\|_{L^{2}\cap W}^2\)\|u(\tau)-v(\tau)\|_{L^{2}\cap W}\;\mathrm{d}\tau.
\end{aligned}
\]
We deduce that
\[
\|\Phi(u)-\Phi(v)\|_{L^{\infty}([0,T];L^{2}\cap W)} \le
C\|\chap{K}\|_{L^{p}}\|u_{0}\|^{2}_{L^{2}\cap W}T\|u-v\|_{L^{\infty}([0,T];L^{2}\cap W)}.
\]
For $T$ possibly smaller (still depending on $\|u_{0}\|_{L^{2}\cap W}$) $\Phi$ is a contraction from $E_{T}$ to $E_{T}$,
so admits a unique fixed point in $E_{T}$, which is a solution for the Cauchy problem.
By resuming the same arguments as in Proposition \ref{prop:uniqueness}
we deduce that the fixed point obtained before is the unique solution
for the Cauchy problem. 

\bigbreak

For $p=1$, the solution constructed before is global in time. In view
of the conservation of the $L^2$-norm, we have
\[
\begin{aligned}
\|u(t)\|_{W}&\le \|u_{0}\|_{W} +
\int_{0}^{t} \|K*|u(\tau)|^{2}\|_{W}\|u(\tau)\|_{L^{2}\cap W}\;\mathrm{d}\tau\\
&\le \|u_{0}\|_{W}+ \int_{0}^{t} C\|\chap{K}\|_{L^{1}}\|u(\tau)\|_{L^{2}}^{2} \|u(\tau)\|_{W}\;\mathrm{d}\tau\\
&\le \|u_{0}\|_{W} + C\|\chap{K}\|_{L^{1}}\|u_{0}\|_{L^{2}}^{2}\int_{0}^{t} \|u(\tau)\|_{W}\;\mathrm{d}\tau,\\
\end{aligned}
\]
and by the Gronwall lemma, we conclude that $\|u(t)\|_{W}$ remains bounded
on finite time intervals. This completes the proof of the first point in
Theorem~\ref{theo:Lp}.

\subsection{Well-posedness in $W$}
Let $\chap{K}\in L^{\infty}$ and consider the Cauchy problem \eqref{eq:r3}
with an initial data in $W$. For $T>0$ we define the following space
\[
Y_{T}=\{u\in C([0,T]; W), \|u\|_{L^{\infty}([0,T];W)} \le 2\|u_{0}\|_{W}\}.
\]
This later is a complete metric space when equipped with the metric
\[
d(u,v)=\|u-v\|_{L^{\infty}([0,T];W)}.
\]
As previously, the local existence of a solution is easily shown by a fixed point argument,
since $\Phi$ is a contraction from $Y_{T}$ to $Y_{T}$, and we show that it is unique.
The proof relies on the following lemma:
\begin{lemma}\label{lem:2}
There exists $C$ such that for all $f,g \in W$,
\[
\|K*(fg)\|_{W}\le C\|\chap{K}\|_{L^{\infty}}\|f\|_{W}\|g\|_{W}.
\]
\end{lemma}

\begin{proof}
Let $f,g\in W$. We have
\[
\begin{aligned}
\|K*(fg)\|_{W} = (2\pi)^{d/2} \|\chap{K}\;\chap{fg}\|_{L^{1}}
&\le (2\pi)^{d/2}\|\chap{K}\|_{L^{\infty}}\|\chap{f}*\chap{g}\|_{L^{1}}\\
&\le
(2\pi)^{d/2}\|\chap{K}\|_{L^{\infty}}\|\chap{f}\|_{L^{1}}\|\chap{g}\|_{L^{1}}\\
&\le
 (2\pi)^{d/2}\|\chap{K}\|_{L^{\infty}}\|f\|_{W}\|g\|_{W}.
\end{aligned}
\]
\end{proof}
It then suffices to reproduce the proof given in the previous
subsection in order to prove the second point of
Theorem~\ref{theo:Lp}, by replacing Lemma~\ref{estWL^2} with
Lemma~\ref{lem:2}.  

\subsection{Ill-posedness in $W$}
For $\gamma \in (0,d)$ we consider the Cauchy problem \eqref{eq:r3}
with the kernel $K$ given by its
 Fourier transform 
 \begin{equation}\label{eq:Khat}
   \chap{K}(\xi)= \frac{1}{|\xi|^{d-\gamma}}{\bf
   1}_{\{|\xi|\le 1\}}.
 \end{equation}
Then, $\chap{K}\in L^{p}$ for all $p \in [1,
\frac{d}{d-\gamma})$. Conversely, for $p\in [1,\infty)$, we can always
find $\gamma \in (0,d)$ such that $K$, defined by \eqref{eq:Khat},
satisfies $\chap{K}\in L^{p}$.

From the first point of Theorem~\ref{theo:Lp}, the Cauchy problem
\eqref{eq:r3} is locally well-posed in $L^2\cap W$.
We suppose that \eqref{eq:r3} is well posed in $W$: let $T>0$ be the local time
existence of the solution. Arguing as in Section~\ref{sec:ill}, this
implies that there exists $C>0$ such that 
for all $f\in W$ and for all $t\in [0, T]$,
\begin{equation}\label{eq:ill2}
\|D(f)(t)\|_{W} \le C\|f\|_{W}^{3},
\end{equation}
where $D$ is the operator defined in Section~\ref{sec:ill}.
Let $f\in \mathcal{S}(\R^d)$. As in Section~\ref{sec:ill} we define the
family of functions 
$(f^{h})_{0<h\le 1}$ by
\[
f^{h}(x) = f(hx).
\]
From \eqref{eq:ill2} we also have for all $t\in [0,T]$
\begin{equation}
\|D(f^{h})(t)\|_{W}\le C\|f^{h}\|^{3}_{W}=C \|f\|^{3}_{W}.
\end{equation}
We define $K_{h}$ by setting
$\chap{K_{h}}(\xi)=\frac{1}{|\xi|^{d-\gamma}} {\bf
  1}_{\{|\xi|>\frac{1}{h}\}}$. We use the following identity:
\[
\begin{aligned}
&\mathcal{F}(K*|U(\tau)f^{h}|^{2}U(\tau)f^{h})(\xi)=\phantom{\int_{\R^{d}}
\int_{\R^{d}}\chap{K}(\xi-y) e^{i\tau|\xi-y-z|^{2}}\chap{f^{h}}(\xi-y-z)
e^{i\tau|y|^{2}}\chap{f^{h}}(y)\;\mathrm{d}y}\\
&=\frac{(2\pi)^{d/2}}{h^{3d}}\iint_{\R^{2d}}
e^{i\tau|\xi-y|^{2}}
e^{i\tau|y-z|^{2}}e^{-i\tau|z|^{2}}\chap{K}(y)
\chap{f}\(\frac{\xi-y}{h}\)\chap{\overline{f}}\(\frac{z}{h}\)
\chap{f}\(\frac{y-z}{h}\)\;\mathrm{d}y\mathrm{d}z\\
&=\frac{(2\pi)^{d/2}}{h^{3d}}\iint_{\R^{2d}}
e^{i\tau|\xi-y|^{2}}
e^{i\tau|y-z|^{2}}e^{-i\tau|z|^{2}}\frac{1}{|y|^{d-\gamma}}
\chap{f}\(\frac{\xi-y}{h}\)\chap{\overline{f}}\(\frac{z}{h}\)
\chap{f}\(\frac{y-z}{h}\)\;\mathrm{d}y\mathrm{d}z\\
&-\frac{(2\pi)^{d/2}}{h^{3d}}\iint_{\R^{2d}}
e^{i\tau|\xi-y|^{2}}
e^{i\tau|y-z|^{2}}e^{-i\tau|z|^{2}}\chap{K_{h}}(y)
\chap{f}\(\frac{\xi-y}{h}\)\chap{\overline{f}}\(\frac{z}{h}\)
\chap{f}\(\frac{y-z}{h}\)\;\mathrm{d}y\mathrm{d}z.\\
\end{aligned}
\]
We inject the above formula into the expression of
$\|D(f^{h})(t)\|_{W}$ and perform  the same changes of variables
as in Section~\ref{sec:ill}. We obtain
\begin{equation*}
\|D(f^{h})(t)\|_{W} \ge \frac{t}{h^{d-\gamma}}\(\|(K\ast|f|^2)f\|_{W} +
 \O\(t h^{2}\)\)- \frac{1}{h^{d-\gamma+2}}X^{h},
\end{equation*}
where 
\[
X^{h}= \left\lVert\int_{0}^{th^{2}}
  U(t-\tau)(K_{h}*|U(\tau)f|^{2}U(\tau)f)\;\mathrm{d}\tau\right\rVert_W.
\]
For $q\in (\frac{d}{d-\gamma}, \infty)$,
\[
\begin{aligned}
X^{h} &\le \int_{0}^{th^{2}}
\left\|\(K_{h}*|U(\tau)f|^{2}\)U(\tau)f\right\|_{W}\;\mathrm{d}\tau 
\le \int_{0}^{th^{2}} \|K_{h}*|U(\tau)f|^{2}\|_{W} \|U(\tau)f\|_{W}\;\mathrm{d}\tau\\
&\le \|f\|_W \int_{0}^{th^{2}}
\|\chap{K_{h}}\|_{L^{q}}
\left\|\F\(|U(\tau)f|^{2}\)\right\|_{L^{q'}}\;\mathrm{d}\tau
\le C th^{2} \|\chap{K_{h}}\|_{L^{q}} ,
\end{aligned}
\]
for some constant $C$ independent of $h\in (0,1]$ and $t\in [0,T]$
(recall that $f\in \Sch(\R^d)$). 
Moreover,
\[
\|\chap{K_{h}}\|_{L^{q}} = \bigg(\int_{\R^{d}}
\frac{1}{|y|^{(d-\gamma)q}}{\bf 1}_{\{|y|>\frac{1}{h}\}}\;\mathrm{d}y\bigg)^{1/q}
= C \bigg(\int_{1/h}^{+\infty}
\frac{r^{d-1}}{r^{(d-\gamma)q}}\;\mathrm{d}r\bigg)^{1/q}, 
\]
so, $\|\chap{K_{h}}\|_{L^{p}}=C h^{d-\gamma-\frac{d}{q}}$ and
\[
\begin{aligned}
\|D(f^{h})(t)\|_{W} & \ge \frac{t}{h^{d-\gamma}}\(\|(K\ast|f|^2)f\|_{W} +
\O\(t h^{2}\)- C h^{d-\gamma-\frac{d}{q}}\).\\
\end{aligned}
\]
Fix $t>0$:
\[
\lim\limits_{h\to +\infty} \frac{t}{h^{d-\gamma}}\(\|(K\ast|f|^2)f\|_{W} +
\O\(t h^{2}\)- C h^{d-\gamma-\frac{d}{q}}\) = +\infty.
\]
So, for $h>0$ small enough we have
\[
\|D(f^{h})(t)\|_{W} > C\|f\|_{W}^{3},
\]
hence a contradiction.

\providecommand{\bysame}{\leavevmode\hbox to3em{\hrulefill}\thinspace}
\providecommand{\MR}{\relax\ifhmode\unskip\space\fi MR }
\providecommand{\MRhref}[2]{%
  \href{http://www.ams.org/mathscinet-getitem?mr=#1}{#2}
}
\providecommand{\href}[2]{#2}

 \end{document}